\documentclass{amsart}
\usepackage{amscd}
\usepackage{amssymb}

\def\Kweb#1{   
http:\linebreak[3]//www.\linebreak[3]math.\linebreak[3]uiuc.
\linebreak[3]edu/\linebreak[3]{K-theory/#1/}}

\def\EE{\mathcal E}
\def\cF{\mathcal F}
\def\cG{\mathcal G}

\def\cO{\mathcal O}
\def\cP{\mathcal P}
\def\cK{\mathcal K}

\def\frakm{\mathfrak m}

\newcommand{\A}{\mathbb{A}}

\newcommand{\Q}{\mathbb{Q}}

\newcommand{\Z}{\mathbb{Z}}

\def\bu{\bullet}

\def\Sing{\operatorname{Sing}}
\def\Spec{\operatorname{Spec}}

\newcommand{\SchF}{\mathrm{Sch}/F}

\newcommand{\SmF}{\mathrm{Sm}/F}

\def\bC{\mathbf C}
\newcommand{\bbH}{\mathbb H}
\newcommand{\bbHcdh}{\bbH_{cdh}}  

\input xypic
\xyoption{all} 

\numberwithin{equation}{section}

\theoremstyle{plain}
\newtheorem{thm}[equation]{Theorem}
\newtheorem{cor}[equation]{Corollary}
\newtheorem{lem}[equation]{Lemma}
\newtheorem{prop}[equation]{Proposition}

\theoremstyle{definition}
\newtheorem{defn}[equation]{Definition}
\theoremstyle{remark}

\newtheorem{ex}[equation]{Example}

\begin{document}
\bibliographystyle{plain}

\title
[$K$-regularity, $cdh$-fibrant $HH$ and a conjecture of Vorst]
{$K$-regularity, \\
$cdh$-fibrant Hochschild homology, \\
and a conjecture of Vorst}

\author{G. Corti\~nas}
\thanks{Corti\~nas' research was partially supported by grants ANPCyT PICT 03-12330, UBACyT-X294, JCyL VA091A05, and MEC MTM00958.}
\address{Dep. Matem\'atica\\ Ciudad Universitaria Pab 1\\ 1428 Buenos Aires, Argentina\\ and  Dep. \'Algebra\\ Fac. de Ciencias\\
Prado de la Magdalena s/n\\ 47005 Valladolid, Spain.}
\email{gcorti@agt.uva.es}\urladdr{http://mate.dm.uba.ar/\~{}gcorti}

\author{C. Haesemeyer}
\address{Dept.\ of Mathematics, University of Illinois, Urbana, IL
61801, USA} \email{chh@math.uiuc.edu}

\author{C. Weibel}
\thanks{Weibel's research was partially supported by NSA grant MSPF-04G-184.}
\address{Dept.\ of Mathematics, Rutgers University, New Brunswick,
NJ 08901, USA} \email{weibel@math.rutgers.edu}

\date{\today}

\begin{abstract}
In this paper we prove that for an affine scheme essentially of
finite type over a field $F$ and of dimension $d$,
$K_{d+1}$-regularity implies regularity, assuming that the
characteristic of $F$ is zero. This verifies a conjecture of
Vorst.
\end{abstract}


\maketitle

\section*{Introduction}

It is a well-known fact that algebraic $K$-theory is homotopy
invariant as a functor on regular schemes; if $X$ is a regular
scheme then the natural map $K_n(X) \to K_n(X\times\A^1)$ is an
isomorphism for all $n\in \Z$. This is false in general for
nonregular schemes and rings.

To express this failure, Bass introduced the terminology that, for
any contravariant functor $\cP$ defined on schemes, a scheme $X$
is called {\it $\cP$-regular} if the pullback maps $\cP(X)\to
\cP(X\times\A^r)$ are isomorphisms for all $r\geq 0$. If $X =
\mathrm{Spec}(R)$, we also say that $R$ is $\cP$-regular. Thus
regular schemes are $K_n$-regular for every $n$. In contrast, it
was observed as long ago as \cite{BM} that a nonreduced affine
scheme can never be $K_1$-regular. In particular, if $A$ is an
Artinian ring (that is, a $0$-dimensional Noetherian ring) then
$A$ is regular (that is, reduced) if and only if $A$ is
$K_1$-regular. In \cite{Vorst2}, Vorst conjectured that for an
affine scheme $X$\!, of finite type over a field $F$ and of
dimension $d$, regularity and $K_{d+1}$-regularity are equivalent;
Vorst proved this conjecture for $d=1$ (by proving that
$K_2$-regularity implies normality).

In this paper, we prove Vorst's conjecture in all dimensions
provided the characteristic of the ground field $F$ is zero. In
fact we prove a stronger statement. We say that $X$ is {\it
regular in codimension $<n$} if $\Sing(X)$ has codimension $\ge n$
in $X$.

Let $\cF_K$ denote the presheaf of spectra such that $\cF_K(X)$ is
the homotopy fiber of the natural map $K(X) \to KH(X)$, where
$K(X)$ is the algebraic $K$-theory spectrum of $X$ and $KH(X)$ is
the homotopy $K$-theory of $X$ defined in \cite{WeibelKH}. We
write $\cF_K(R)$ for $\cF_K(\mathrm{Spec}(R)).$

\begin{thm}\label{thm:main-intro}

Let $R$ be a commutative ring which is essentially of finite type
over a field $F$ of characteristic $0$. Then:
\newline (a) If
$\cF_K(R)$ is  $n$-connected, then $R$ is regular in codimension
$<n$.
\newline
(b) If $R$ is $K_n$-regular, then $R$ is regular in codimension $<n$.
\newline
(c) (Vorst's conjecture) If $R$ is $K_{1+\dim(R)}$-regular,
then $R$ is regular.
\end{thm}

It was observed in \cite{WeibelKH} that if $X$ is $K_n$-regular
then $K_i(X)\to KH_i(X)$ is an isomorphism for $i\le n$, and a
surjection for $i=n+1$, so that $\cF_K(X)$ is $n$-connected. Thus
(a) implies (b) in this theorem, and (c) is a special case of (b).

The bounds in (a) and (b) are the best possible, because it
follows from Vorst's results (\cite[Thm. A]{Vorst2}, \cite[Thm.
3.6]{Vorst1}) that for an affine singular seminormal curve $X$,
$\cF_K(X)$ is $1$-connected, but not $2$-connected. The converse
of (c) is trivial, but those of (a) and (b) are false. Indeed,
affine normal surfaces are regular in codimension~1 but may not be
$K_{-1}$-regular, much less $K_2$-regular; see
\cite[5.8.1]{WeibelNorm}.

Finally the analogue of (c) --and thus also of (a) and (b)-- for
$K$-theory of general nonaffine schemes is false. Indeed we give
the following example of a nonreduced (and in particular
nonregular) projective curve which is $K_n$-regular for all $n.$

\begin{thm}\label{thm:counteregg-intro}
Let $(X,Q)$ be an elliptic curve over a field of characteristic~0,
and $P$ a rational point on $X$ such that the line bundle
$L=L(P-Q)$ does not have odd order in the Picard group $Pic(X)$.
Write $Y$ for the nonreduced scheme with the same underlying space
as $X$ but with structure sheaf $\cO_Y=\cO\oplus L$, where $L$ is
regarded as a square-zero ideal.

\noindent Then $Y$ is $K_n$-regular for all $n$, and $\cF_K(Y)$ is
contractible.
\end{thm}

The proof of Theorem \ref{thm:main-intro} employs results from our
paper with M. Schlichting \cite{chsw} that allow us to describe
$\cF_K$ in terms of cyclic homology; the necessary statements will
be recalled in Section \ref{sec:recoll}. In Section
\ref{sec:cdhHH}, we study the $cdh$-fibrant version of Hochschild
homology and its Hodge decomposition. Section \ref{sec:HHsmooth}
contains a smoothness criterion (Theorem \ref{thm:HHcriterion})
using $cdh$-fibrant Hochschild homology, which is of independent
interest and generalized in Theorem \ref{thm:HHcodim-h}. The proof
of part (a) of Theorem \ref{thm:main-intro} is contained in
Section \ref{sec:main} (see Theorem \ref{thm:main}). As explained
above, parts (b) and (c) follow from this. Finally Section
\ref{sec:counteregg} is devoted to the counterexample stated in
Theorem \ref{thm:counteregg-intro} (and restated as Theorem
\ref{thm:counteregg}).

\medskip
{\bf Notation}

All rings considered in this paper are commutative and noetherian.
We shall write $\SchF$ for the category of schemes essentially of finite
type over a field $F$. Objects of $\SchF$ shall be called {\it $F$-schemes}.

The category of spectra we use in this paper will not be critical.
In order to minimize technical issues, we will use the terminology
that a {\it spectrum} $\EE$ is a sequence of simplicial sets
$\EE_n$ together with bonding maps $b_n: \EE_n \to
\Omega\EE_{n+1}$. We say that $\EE$ is an {\it $\Omega$-spectrum}
if all bonding maps are weak equivalences. A map of spectra is a
strict map. We will use the model structure on the category of
spectra defined in \cite{BF}. Note that in this model structure,
every fibrant spectrum is an $\Omega$-spectrum.

We say that a presheaf $E$ of spectra on $\SchF$ satisfies the
{\it Mayer-Vietoris-property} (or MV-property, for short) for a
cartesian square of schemes
\begin{equation}\label{square}
\begin{CD}
Y' @>>> X' \\
@VVV @VVV \\
Y @>>> X
\end{CD}
\tag{$\square$}
\end{equation}
if applying $E$ to this square results in a homotopy cartesian
square of spectra.  We say that $E$ satisfies the Mayer-Vietoris
property for a class of squares provided it satisfies the
MV-property for each square in the class.

We say that $E$ satisfies {\it Nisnevich descent} for $\SchF$ if
$E$ satisfies the MV-property for all elementary Nisnevich squares
in $\SchF$; an {\em elementary Nisnevich square} is a cartesian
square of schemes \eqref{square} for which $Y\to X$ is an open
embedding, $X'\to X$ is \'etale and $(X'-Y')\to(X-Y)$ is an
isomorphism.  By \cite[4.4]{Nis}, this is equivalent to the
assertion that $E(X) \to \bbH_{nis}(X,E)$ is a weak equivalence
for each scheme $X$, where $\bbH_{nis}(-,E)$ is a fibrant
replacement for $E$ in a suitable model structure.

We say that $E$ satisfies {\it $cdh$-descent} for $\SchF$ if $E$
satisfies the MV-property for all elementary Nisnevich squares
(Nisnevich descent) and for all abstract blow-up squares in
$\SchF$; an {\em abstract blow-up square} is a square
\eqref{square} such that $Y\to X$ is a closed embedding, $X'\to X$
is proper and the induced morphism $(X' - Y')_{red} \to
(X-Y)_{red}$ is an isomorphism. With M. Schlichting, we showed in
Theorem~3.7 of \cite{chsw} that $cdh$-descent is equivalent to the
assertion that $E(X) \to \bbHcdh(X,E)$ is a weak equivalence for
each scheme $X$, where $\bbHcdh(-,E)$ is a fibrant replacement for
the presheaf $E$ in a suitable model structure. We abbreviate
$\bbHcdh(-E)$ as $\bbH(-,E)$ if no confusion can arise, and write
$\bbH_n(X,E)=\bbH^{-n}(X,E)$ for $\pi_n \bbH(X,E)$.

We use cohomological indexing for all chain complexes in this
paper; for a complex $A$, $A[p]^q = A^{p+q}$. It is well-known
that there is an Eilenberg-Mac{\,}Lane functor $A\mapsto|A|$ from
chain complexes of abelian groups to spectra, and from presheaves
of chain complexes of abelian groups to presheaves of spectra.
This functor sends quasi-isomorphisms of complexes to weak
homotopy equivalences of spectra, and satisfies
$\pi_n(|A|)=H^{-n}(A)$. In this spirit, we will use descent
terminology for presheaves of complexes.

For example, the Hochschild, cyclic, periodic and negative cyclic
homology of schemes over a field $k$ (such as $F$-schemes over a
field $F\supseteq k$) can be defined using the Zariski
hypercohomology of certain presheaves of complexes; see
\cite{WeibelHC} and \cite[2.7]{chsw} for precise definitions.  We
shall write these presheaves as $HH(/k)$, $HC(/k)$, $HP(/k)$ and
$HN(/k)$, respectively and regard them as presheaves of either
cochain complexes or spectra. When $k$ is omitted, it is understood
that $k=\Q$ is intended.
Finally, we write $\Omega^i_{/k}$
for the presheaf $X\mapsto \Omega^i_{X/k}$, while
$\Omega^\bu_{/k}$ denotes the presheaf of algebraic de Rham
complexes and $\Omega^{\leq i}_{/k}$ denotes its brutal truncation
in degree $i$.

\section{$cdh$-descent}\label{sec:recoll}

In this section we recall the main results from \cite{chsw}, and
prove that the failure of $K$-theory to be homotopy invariant can
be measured using cyclic homology. We work on the category $\SchF$
of $F$-schemes essentially of finite type over a field $F$ of
characteristic $0$.

Here are two of the main results of \cite{chsw}. Recall that
infinitesimal $K$-theory $\cK^{inf}(X)$ is the homotopy fiber of
the Jones-Goodwillie Chern character $K(X) \to HN(X)$. The first
one is Theorem~4.6 of \cite{chsw}:

\begin{thm}\label{thm:Kinfdesc}
The presheaf of spectra $\cK^{inf}$ satisfies $cdh$-descent.
\end{thm}

The second one is a slight modification of Corollary~3.13 of
\cite{chsw}.

\begin{thm}\label{thm:HPdesc}
For each subfield $k \subseteq F$, the presheaf
$HP(/k)$ satisfies $cdh$-descent on $\SchF$.
In particular, $HP$ satisfies $cdh$-descent.
\end{thm}

\begin{proof}
For $k=F$, this is proven in \cite[Corollary 3.13]{chsw}.
As in {\it loc.\ cit.,} the result for $k\subset F$ follows from
\cite[Theorem 3.12]{chsw}, once we verify that the hypotheses hold.
But this follows from three observations:
(1) the Cuntz-Quillen excision theorem holds over $k$
(see \cite[5.3]{CQ}, noting that the condition that the base field
be algebraically closed is not needed, see \cite[p.\ 3]{CQ});
(2) Goodwillie's theorem that periodic cyclic homology is invariant
under infinitesimal extension does not require rings of finite type
over $k$ (see \cite[II.5.1]{Goodw} or \cite[4.1.15]{LodayHC92}); and
(3) the results of \cite[Section 2]{chsw} also hold for
(Hochschild, cyclic, periodic, negative) homology over $k$.
\end{proof}

\begin{ex}\label{ex:KH-cdh}
Consider the presheaf of spectra $KH$ associated to Homotopy
$K$-theory. The main theorem of \cite{HKH} states that
$\bbHcdh(-,K)\simeq KH$, and that $KH$ satisfies $cdh$ descent.
Thus the following definition of $\cF_K$ is compatible with the
definition of $\cF_K$ given in the introduction.
\end{ex}

\begin{defn}\label{def:ndescent}
For any presheaf of spectra $E$, we write $\cF_E$ for the homotopy
fiber of $E\to \bbHcdh(-,E)$. If $\cF_E(X)$ is $n$-connected
for some scheme $X$ (for all $X$ in some subcategory of $Sch/F$),
we say that $E$ satisfies $n$-$cdh$-descent on $X$ (resp., on the
subcategory). Note that if $E$ satisfies $n$-$cdh$-descent
for all $n$ on a subcategory, then $E$ satisfies $cdh$-descent on
that subcategory.
\end{defn}

\smallskip\noindent
Since the fibrant replacement functor $\bbHcdh$ preserves
(objectwise) fibration sequences, it follows that $\cF$ does too. (See the
first paragraph of \cite[Section 5]{chsw}.) We record this as the
following observation.

\begin{lem}\label{lem:3by3}
Let $E_1\to E_2\to E_3$ be a fibration sequence of presheaves of
spectra. Then there is a natural induced fibration sequence
\[ \cF_{E_1} \to \cF_{E_2} \to \cF_{E_3}. \]
\end{lem}

\begin{thm}\label{thm:fibers3}
For any scheme $X$, essentially of finite type over a field of
characteristic $0$, the Chern character $K\to HN=HN(/\Q)$ induces
natural weak equivalences
$$\cF_K (X) \xrightarrow{\simeq} \cF_{HN}(X)
\xleftarrow{\simeq} \Omega^{-1}\cF_{HC}(X).$$
\end{thm}

\begin{proof} The first weak equivalence follows from  Lemma \ref{lem:3by3}
and Theorem \ref{thm:Kinfdesc}. The second weak equivalence
follows from  Lemma \ref{lem:3by3}, Theorem \ref{thm:HPdesc} and
the $SBI$ fibration sequence $\Omega HP \to \Omega^{-1}HC \to HN\to HP$.
\end{proof}

\begin{cor}\label{lem:main1}
Let $X\in\SchF$. Then $K$ satisfies $(n+1)$-$cdh$-descent on $X$
if and only if $HH$ satisfies $n$-$cdh$-descent on $X$.
\end{cor}

\begin{proof}
Since  $HH(X)$ and $\bbH(X,HH)$ are
$n$-connected for $n<-\dim(X)$ by \cite{WeibelHC},
this follows from Theorem \ref{thm:fibers3} and the SBI sequence
$\Omega^{-1}HC \to HH \to HC$.
\end{proof}

\bigskip
\section{$cdh$-fibrant Hochschild and cyclic homology.}\label{sec:cdhHH}

In this section we study the $cdh$-fibrant version of Hochschild
homology and its Hodge decomposition, and establish some of their
basic properties.

For legibility, we will write $a$ for the natural morphism of
sites from the $cdh$-site to the Zariski site on $\SchF$.  If $A$
is a Zariski sheaf, its associated $cdh$ sheaf will be written as
$A_{cdh}$ or $a^* A$. We will simplify notation and write
$H^*_{cdh}(X,A)$ for the cohomology of $A_{cdh}$. In particular
this applies to the sheaf $X\mapsto\Omega^i_{X/k}$ of K\"ahler
$i$-differential forms ($i\ge 0$); $H_{cdh}^*(X,\Omega^i_{/k})$ is
the cohomology of $a^*\Omega^i_{/k}$. If $A$ is a complex of
presheaves of abelian groups on $\SchF$, then we write
$\bbHcdh(A)$ for a $cdh$-fibrant replacement of $A$, and
$\bbH_{cdh}(X,A)$ for its complex of sections over $X$; the usual
hypercohomology $\bbHcdh^n(X,A)$ is just $H^{n}\bbHcdh(X,A)$. For
example, if $A$ is a presheaf, considered as a complex
concentrated in degree zero, then $\bbHcdh(A)$ is just an
injective resolution of the $cdh$-sheafification $A_{cdh}$, and
$\bbHcdh^n(X,A)$ is the usual cohomology $H_{cdh}^n(X,A_{cdh})$ of
$A_{cdh}$.

When $A$ is an unbounded complex, such as a complex representing
Hochschild homology, then $\bbH_{cdh}(A)$ may be constructed using
product total complexes of flasque Cartan-Eilenberg resolutions.
This works because the columns of the Cartan-Eilenberg double
complex are locally cohomologically bounded by \cite{SVBK}.

The $cdh$ site is Noetherian (every covering has a finite
subcovering), so $cdh$ cohomology commutes with filtered direct
limits of sheaves. A typical application of this fact is that if
$M$ is a sheaf of $F$-modules and $V$ is a vector space then
$H_{cdh}^n(X,V\otimes_F M)\cong V\otimes_F H_{cdh}^n(X,M)$.  (See
\cite{SGA4II}, Exp.~VI, 2.11 and 5.2.)

The Hochschild and cyclic homology of schemes over a field $k$
(such as $F$-schemes over a field $F\supseteq k$) can be defined
using the Zariski hypercohomology of certain presheaves of mixed
complexes; see \cite{WeibelHC} and \cite[2.7]{chsw} for precise
definitions. It was observed in \cite[3.0]{WeibelHo} that, because
the mixed complexes already admit a Hodge decomposition, so do the
complexes $HH(/k)$, $HC(/k)$, $HP(/k)$ and $HN(/k)$. Taking
fibrant replacements for any Grothendieck topology respects such
product decompositions; the following proposition records this for
the $cdh$-topology.

\begin{prop}\label{prop:cdhhodge}
Let $X$ be a scheme over a field $k$ of characteristic $0$ (such as an
$F$-scheme for $F\supseteq k$).
Then the $cdh$-fibrant Hochschild, cyclic, negative cyclic and
periodic cyclic homology of $X$ over $k$ admit natural Hodge
decompositions.

\noindent That is, if $H$ denotes any of $HH(/k)$, $HC(/k)$,
$HP(/k)$ or $HN(/k)$ then:
$$\bbHcdh(X,H) \cong \prod\bbHcdh(X,H^{(i)}).$$
\end{prop}

Moreover, using the computations of these decompositions provided
in \cite[Theorem 3.3]{WeibelHo} and the fact that all $F$-schemes
are locally smooth in the $cdh$-topology, it is possible to
compute the Hodge decomposition explicitly in terms of the
$cdh$-hypercohomology of the de Rham complex.

\begin{thm}\label{thm:cdhhodge}
Let $k\subseteq F$ be a subfield. There are natural isomorphisms for
every $F$-scheme $X$:
\begin{equation*}
\begin{split}
\pi_n\bbHcdh(X,HH^{(i)}(/k)) \cong H^{i-n}_{cdh}(X,\Omega^i_{/k});\\
\pi_n\bbHcdh(X,HC^{(i)}(/k)) \cong \bbHcdh^{2i-n}(X,\Omega_{/k}^{\leq i});\\
\pi_n\bbHcdh(X,HN^{(i)}(/k)) \cong \bbHcdh^{2i-n}(X,\Omega_{/k}^{\geq i});\\
\pi_n\bbHcdh(X,HP^{(i)}(/k)) \cong \bbHcdh^{2i-n}(X,\Omega_{/k}^\bu).
\end{split}
\end{equation*}
\end{thm}

\begin{proof}
Let $\bC(X)$ denote the mixed complex computing the Hochschild and
cyclic homology of $X$ over $k$. The functor $\bC:X\mapsto \bC(X)$
is a presheaf of mixed complexes. By \cite[9.8.12]{WeibelHA94},
there is a Hodge decomposition $\bC \cong \prod_i \bC^{(i)}$ and a
natural map of mixed complexes $e:\bC\to (\Omega^*_{X/k},0,d)$
that sends the Hochschild chain complex
$HH^{(i)}(X)=(\bC^{(i)}(X),b)$ to $\Omega^i_{X/k}[i]$.  As
observed in \cite{WeibelHo}, the induced map on Connes' double
complexes sends $\mathcal{B}_{**}^{(i)}$ to $\Omega_{X/k}^{\leq
i}[2i]$. It suffices to prove that these are locally
quasi-isomorphisms for the $cdh$ topology. This boils down to
showing that $e$ induces a quasi-isomorphism
$HH^{(i)}(R/k)\to\Omega^i_{R/k}[i]$ for every regular local
$F$-algebra $R$. For $k=F$, this is the
Hochschild-Kostant-Rosenberg theorem (\cite[9.4.7]{WeibelHA94}).
The general case follows from this and the fact that $R$ is the
union of smooth $k$-algebras. (It also follows from the
Kassel-Sledsj\oe\ spectral sequence of \cite[4.3a]{kasle}, which
we recall in \ref{lem:main11} below.)
\end{proof}

\goodbreak
\begin{lem}\label{lem:tech1}
Let $R$ be an $F$-algebra essentially of finite type, $k\subseteq
F$ a subfield. Then for each $n$, $HH(/k)$ satisfies $n$-cdh
descent on $X=\Spec(R)$ if and only if the following three
conditions hold simultaneously:
\addtocounter{equation}{-1}
\begin{subequations}
\begin{gather}
HH_m^{(q)}(R/k)=0 \text{\quad if } 0\le q<m\le n;\label{hhkvanish}\\
\Omega^q_{R/k}\rightarrow H^0_{cdh}(X,\Omega^q_{/k})\quad
\text{is bijective if } q\le n \text{ and onto if } q=n+1;\label{agreek}\\
H^p_{cdh}(X,\Omega^q_{/k})=0 \quad\text{ if } p>0 \text{ and }0\le q\le p+n+1.
\label{cdhkvanish}
\end{gather}
\end{subequations}
\end{lem}

\noindent Note that \eqref{hhkvanish} is vacuous if $n\le0$, and
\eqref{agreek} is vacuous if $n\le-2$. In particular, $HH(/k)$
satisfies $(-2)$-cdh-descent just in case
$H^p_{cdh}(X,\Omega^q_{/k})=0$ for all $p>q\ge0$.

\begin{proof} This follows easily from the Hodge decomposition
and the isomorphisms
\[HH_q^{(q)}(R/k) \cong \Omega^q_{R/k},
\quad\hbox{and}\quad HH_m^{(q)}(R/k) = 0 \hbox{\ for } q>m.\] In
more detail, we see from \ref{prop:cdhhodge} and
\ref{thm:cdhhodge} that the maps $HH_m^{(q)}(R/k)\to
H_{cdh}^{q-m}(X,\Omega^q_{/k})$ must be isomorphisms for $m\le n$
and onto for $m=n+1$.
\end{proof}

On smooth schemes, all our functors are well-behaved. Recall from
\cite{SVBK} that the $scdh$ topology on $\SmF$ is the restriction
of the $cdh$ topology on $\SchF$. Since every scheme is locally
smooth, it follows that $\bbH_{scdh}(X,A)$ is just $\bbHcdh(X,A)$
for every presheaf $A$. (See the argument of the first part of the
proof of \cite[3.12]{chsw}.) If $A$ satisfies $scdh$-descent then
$A(X)\cong\bbHcdh(X,A)$ for all smooth $X$.

Recall that $k$ is a subfield of $F.$

\begin{thm}\label{thm:smoothcdh}
Let $H$ denote any of: $HH(/k)$, $HC(/k)$, $HN(/k)$ or $HP(/k)$,
and let $H^{(i)}$ denote the $i^{th}$ component in the Hodge
decomposition of $H$.

\noindent Then $H$ and $H^{(i)}$ satisfy $scdh$-descent on $\SmF$.
In particular, if $X$ is smooth over $F$, then $H^{(i)}(X) \cong
\bbHcdh(X,H^{(i)})$.
\end{thm}

\begin{proof} Since every smooth scheme over $F$ is locally a filtered limit
of smooth affine schemes over $k$, and  $H$ commutes with limits
of affine schemes, we may assume that $k=F$. By \cite [3.9, 2.9,
and 2.10]{chsw}, Hochschild, cyclic, negative and periodic cyclic
homology (relative to $F$) all satisfy $scdh$-descent on $\SmF$.

By Proposition~\ref{prop:cdhhodge}, the quasi-isomorphisms
$H(X)\cong\bbHcdh(X,H)=\bbH_{scdh}(X,H)$ induce quasi-isomorphisms
$H^{(i)}(X)\cong \bbH_{cdh}^{(i)}(X,H) = \bbH_{scdh}^{(i)}(X,H)$
for all $i$.
\end{proof}

The special case $H^*_{Zar}(X,\cO) \cong H^*_{cdh}(X,\cO)$ (for smooth $X$)
of the following corollary was proven in \cite[6.3]{chsw}.

\begin{cor}\label{cor:Omega-scdh}
If $X$ is smooth over $F$, then $H^p_{Zar}(X,\Omega^i_{/k}) \cong
H^p_{cdh}(X,\Omega^i_{/k})$ for all $p$ and $i$. In particular,
$\Omega^i_{X/k}\cong H_{cdh}^0(X,\Omega^i_{/k})$.
\end{cor}

\begin{proof} Consider the map
$e^{(i)}:HH^{(i)}=(\bC^{(i)},b) \to \Omega^i_{/F}[i]$
 of complexes of Zariski sheaves.  By \cite[3.3]{WeibelHo}, it
is a quasi-isomorphism over every smooth scheme $X$ over $F$,
inducing $HH_{i-p}^{(i)}(X)\cong H^{p-i}_{Zar}(X,\bC^{(i)})\cong
H^p_{Zar}(X,\Omega^i_{/k})$. The map $e^{(i)}$ remains a
quasi-isomorphism after sheafifying for the $cdh$ topology, so
that $\bbHcdh^{p-i}(X,HH^{(i)})\cong H^{p-i}_{cdh}(X,\bC^{(i)})
\cong H^p_{cdh}(X,\Omega^i_{/k})$. By Theorem \ref{thm:smoothcdh},
$HH_{n}^{(i)}(X)\cong \bbHcdh^{-n}(X,HH^{(i)})$, whence the
result.
\end{proof}

The next result is proven by copying the proof of
\cite[6.1]{chsw}, replacing $\cO$ with $\Omega^i_{/k}$.

\begin{prop}\label{prop:Hdzero}
If $X$ is a $d$-dimensional scheme, essentially of finite
type over $F$, and $k\subseteq F$ is a subfield, then
\[ H^d_{Zar}(X,\Omega^i_{/k}) \to H^d_{cdh}(X,\Omega^i_{/k}) \]
is surjective. In particular, if $X$ is affine and $d>0$ then
$H^d_{cdh}(X,\Omega^i_{/k})=0$.
\end{prop}

The following useful theorem is proven in \cite[12.1]{SVBK}.

\begin{thm}\label{cdh:MV}
For every abstract blow-up square $(\square)$, and for every
complex of sheaves of abelian groups $A$, there is a long exact
Mayer-Vietoris sequence:
\[
\cdots H^n_{cdh}(X,A) \to H^n_{cdh}(X',A)\oplus H^n_{cdh}(Y,A) \to
H^n_{cdh}(Y',A) \to H^{n+1}_{cdh}(X,A) \cdots
\]
\end{thm}

Consider the change-of-topology morphism
$a:(\SchF)_{cdh}\to(\SchF)_{Zar}$.

\begin{lem}\label{lem:qc}
If a Zariski sheaf $M$ on $\SchF$ is a quasi-coherent sheaf
(resp., coherent sheaf) on each $X_{Zar}$, and $M$ satisfies
$scdh$-descent on $\SmF$, then the cohomology sheaves
$R^qa_*(a^*M)$ are also quasi-coherent (resp., coherent) on each
$X_{Zar}$.

If $X=\Spec(R)$ is affine, then $R^qa_*(a^*M)$ is the quasi-coherent
sheaf associated to the $R$-module $H^q_{cdh}(X,M)$, and the natural
map $M(X)\to H^0_{cdh}(X,M)$ is $R$-linear.
\end{lem}

\begin{proof}
We proceed by induction on $\dim X$, the case $\dim(X)=0$ being
clear. Pick a smooth proper birational cdh cover $X'$ of $X$, as
in \cite[5.9]{SVBK} or \cite[12.23]{MVW}, and form the abstract
blow-up square $(\square)$.  By Theorem \ref{cdh:MV}, we get a triangle
on $X_{Zar}$:
$Ra_*(a^*M)|_X \to Ra_*(a^*M)|_{X'\coprod Y} \to Ra_*(a^*M)|_{Y'}$
As the latter two terms have quasi-coherent (resp., coherent) cohomology
sheaves, by induction and $scdh$-descent on $X'$, so does the first.

If $X$ is affine, then $H^p_{Zar}(X,R^qa_*M)=0$ for $p>0$. Hence
the Leray spectral sequence collapses to yield
$H^q_{cdh}(X,M)=H^0_{Zar}(X,R^qa_*(a^*M))$.
\end{proof}

\begin{cor}\label{cor:qc}
Suppose that $X = \mathrm{Spec}(R)$ is affine. Then
$U\mapsto\pi_n\cF_{HH(/k)}(U)$ and
$U\mapsto\pi_n\bbHcdh(U,HH(/k))$ are quasi-coherent Zariski
sheaves on $X$ for all $n$.
\end{cor}

\section{A criterion for smoothness.}\label{sec:HHsmooth}

In this section we present a local criterion for smoothness of
schemes over a field $F$, in terms of the Hochschild homology and
$cdh$-fibrant Hochschild homology of their local rings over $F$
(see \ref{thm:HHcriterion}). As an application we prove Vorst's
conjecture for algebras of finite type over $\Q$ and their
localizations at maximal ideals (see \ref{thm:mainQ}).

A stronger global version of the following result shall be proved
in Section \ref{sec:main} below (see Theorem \ref{thm:HHcodim-h}).

Recall that $F$ is a field of characteristic $0.$

\begin{thm}\label{thm:HHcriterion}
Let $R$ be the local ring of a $d$-dimensional $F$-algebra of
finite type at a maximal ideal. If $HH(/F)$ satisfies
$d$-$cdh$-descent on $R$, then $R$ is smooth over $F$.
\end{thm}

\begin{proof}
Recall that $\Omega_{/F}^\bu$ denotes the de Rham complex, whose
terms are the Zariski sheaves $\Omega^i_{/F}$, while
$\Omega_{/F}^{\leq i}$ denotes its brutal truncation in degrees at
most $i$. By \ref{thm:cdhhodge} and \ref{thm:HPdesc}, we have
isomorphisms
\[\xymatrix{
HP^{(j)}_n(X/F)\ar[r]^(0.4){\cong}&\bbHcdh^{2j-n}(X,\Omega^\bullet_{/F})}
\]
for any $X\in \SchF$, and all $n$ and $j$. Moreover, by the proof
of \ref{thm:cdhhodge}, this isomorphism factors through a natural
map $e:HP^{(j)}_n(X/F)\to
\mathbb{H}_{Zar}^{2j-n}(X,\Omega^\bullet_{/F})$. Now specialize to
the case $X=\Spec R$, where $R$ is as in the theorem. Since every
$cdh$ cover of $X$ has a $d$-dimensional refinement, we have
\[
\bbHcdh^{*}(X,\Omega^\bullet_{/F})=\bbHcdh^{*}(X,\Omega^{\le d}_{/F}).
\]
Moreover, Lemma \ref{lem:tech1} implies that the hypercohomology
spectral sequence for $\bbHcdh^*$ degenerates to yield an isomorphism
\[
H^*(\Omega^{\le d}_{R/F},d)\to \bbHcdh^{*}(X,\Omega^{\le d}_{/F}).
\]
The canonical map $S:HP^{(j)}_n(R/F)\to HC^{(j-1)}_{n-2}(R/F)$ fits into the
commutative diagram
\[ \xymatrix{
HP_{d+2}^{(d+1)}(R/F)\ar[d]_S\ar[r]_e&
H_{dR}^d(R/F)\ar@{^{(}->}[d]\ar[r]&
\bbHcdh^{d}(X,\Omega^{\bullet}_{/F})\ar[d]^{\wr}\\
          HC_d^{(d)}(R/F)\ar[r]^{\cong}_e&
\Omega^d_{R/F}/d\Omega^{d-1}_{R/F}\ar[r]^{\cong}&
\bbHcdh^{d}(X,\Omega^{\le d}_{/F}).}
\]
We have seen that the top composite, the right vertical and both
bottom arrows are isomorphisms. It follows that the middle vertical
inclusion is the identity map, i.e., that $d\Omega^d_{R/F}=0$.
On the other hand, $d\Omega^{d}_{R/F}$ generates $\Omega^{d+1}_{R/F}$ as an
$R$-module; therefore we can infer that $\Omega^{d+1}_{R/F} = 0$.
By Lemma \ref{lem:tech2} below, $R$ is regular, and hence smooth
over $F$.
\end{proof}

\begin{lem}\label{lem:tech2} Let $F$ be any perfect
field. Suppose $R$ is the local ring of a $d$-dimensional
$F$-algebra of finite type at a maximal ideal. If
$\Omega^{d+1}_{R/F} = 0$, then $R$ is regular.
\end{lem}

\begin{proof} Let $\frakm$ be the maximal ideal of $R$.
Since $L:=R/\frakm$ is smooth over $F$, the Second Fundamental
Theorem \cite[9.3.5]{WeibelHA94} shows that there is an
isomorphism $\frakm/\frakm^2\to L\otimes_R\Omega^1_{R/F}$ sending
$x$ to $dx$. Consequently, there is a surjection from
$\Omega^{d+1}_{R/F}$ onto $\Lambda^{d+1}_L(\frakm/\frakm^2)$,
which is a nonzero vector space unless $R$ is regular.
\end{proof}

As an application, we can now verify Vorst's Conjecture for
algebras of finite type over $\Q$ and their localizations at
maximal ideals.

\begin{thm}\label{thm:mainQ}
Let $R$ be a $d$-dimensional $\Q$-algebra which is either of
finite type over $\Q$, or a localization of a $\Q$-algebra of
finite type at a maximal ideal.

If $R$ is $K_{d+1}$-regular, then $R$ is regular.
\end{thm}

\begin{proof}
First assume $R$ is of finite type over $\Q$, and $K_{d+1}$-regular.
To prove $R$ is regular, we may replace $R$ by its localization at a
maximal ideal; these local rings are also $K_{d+1}$-regular,
by Vorst's localization theorem \cite[1.90]{Vorst1}.
Thus we are reduced to proving the theorem in the local case.

As remarked in the introduction, if $R$ is $K_{d+1}$-regular, then
$\cF_K(R)$ is $(d+1)$-connected (see \cite{WeibelKH}).
By Corollary \ref{lem:main1}, $\cF_{HH(/F)}(R)$ is $d$-connected. Now
Theorem \ref{thm:HHcriterion} applies to prove that $R$ is smooth
over $\Q$ and hence regular.
\end{proof}

\section{Vorst's conjecture.}\label{sec:main}

In this section we will prove Theorem \ref{thm:main-intro}. Throughout,
$F$ will be a fixed field of characteristic zero, $k\subseteq F$ a subfield,
$R$ an $F$-algebra essentially of finite type, and $X=\Spec(R)$. Note
that we write $HH$ for $HH(/\Q)$.

\begin{lem}\label{lem:main11} (Kassel-Sledsj\oe, \cite[4.3a]{kasle})
Let $k\subseteq F$ and $p\ge1$ be fixed. Then there is a bounded
second quadrant homological spectral sequence ($0\le i<p$,
$j\ge0$):
\[
{}_p E^1_{-i,i+j}=\Omega^{i}_{F/k}\otimes_FHH^{(p-i)}_{p-i+j}(R/F)
\Rightarrow HH_{p+j}^{(p)}(R/k)
\]
\end{lem}

\begin{lem}\label{lem:main12} Let $k\subseteq F$ and $p\ge1$ be fixed.
Then there is a spectral sequence:
\[
{}^p
E_1^{i,j}=\Omega^i_{F/k}\otimes_FH^{i+j}_{cdh}(X,\Omega^{p-i}_{/F})
\Rightarrow H^{i+j}_{cdh}(X,\Omega^{p}_{/k}).
\]
\end{lem}
\begin{proof}
Consider the sheaf of ideals $I:=\ker(\Omega^*_{/k}\to
\Omega^*_{/F})$. The $I$-adic filtration of $\Omega^*_{/k}$
induces a filtration $\cG=\cG(p)$ on $\Omega^p_{/k}$. If $R$ is
any $F$-algebra essentially of finite type, we have a natural
surjection
\begin{equation}\label{filt}
\Omega^i_{F/k}\otimes_F\Omega^{p-i}_{R/F}
    \twoheadrightarrow \cG^i(R)/\cG^{i+1}(R),
\end{equation}
which is an isomorphism if $R$ is smooth. Thus the
cdh-sheafification of \eqref{filt} is an isomorphism. Since
$\Omega^i_{F/k}$ is a vector space, the spectral sequence of the
lemma is the one associated to the corresponding filtration of the
$cdh$ sheaf $a^*\Omega^p_{/k}$.
\end{proof}

\begin{lem}\label{lem:main13}
Let $X=\Spec(R)$ be affine, and fix $n\ge 0$. Assume that
\addtocounter{equation}{-1}
\begin{subequations}
\begin{gather}
\Omega^q_{R/k}\rightarrow H^0_{cdh}(X,\Omega^q_{/k})\quad
\text{is bijective if } q\le n \text{ and onto if } q=n+1,\label{preagree}\\
H^1_{cdh}(X,\Omega^q_{/F})=0\text{ if } q\le n+1.\label{1cdhFvanish}
\end{gather}
\end{subequations}
Then $\Omega^q_{R/F}\rightarrow H^0_{cdh}(X,\Omega^q_{/F})$ is
bijective if $q\le n$, and onto if $q=n+1$.
\end{lem}

\begin{proof}
By induction on $q$. If $q=0$, there is nothing to prove.
Fix $q>0$, and consider the filtration $\cG^i = \cG^i(q)$, $0\le i\le q$,
considered in the proof of Lemma \ref{lem:main12}.
We have a commutative diagram
\begin{equation}\label{preinfty}
\xymatrix{\Omega^i_{F/k}\otimes_F\Omega^{q-i}_{R/F}\ar[d]\ar[r]&
 {}^{}\cG^i(R)/ \cG^{i+1}(R)\ar[d]\\
         \Omega^i_{F/k}\otimes_FH^0_{cdh}(X,\Omega^{q-i}_{/F})\ar[r]&
  H^0_{cdh}(X,\cG^i/\cG^{i+1}).}
\end{equation}
The top arrow is surjective for all $i$, and an isomorphism for
$i=0$. The bottom arrow is an isomorphism by the proof of
\ref{lem:main12}. By the inductive hypothesis, the left vertical arrow
is an isomorphism for $0<i$. It follows that the top
arrow is an isomorphism for all $0\le i\le q$, and that the arrow
on the right is an isomorphism
for $i>0$. By \eqref{1cdhFvanish} we have an exact sequence:
\[
\xymatrix{0\ar[r]&H^0_{cdh}(X,\cG^{i+1})\ar[r]&H^0_{cdh}(X,\cG^{i})\ar[r]&
\Omega^i_{F/k}\otimes_FH^0_{cdh}(X,\Omega^{q-i}_{/F})\ar[d]\\
     &\qquad\qquad0&H^1_{cdh}(X,\cG^{i})\ar[l]&H^1_{cdh}(X,\cG^{i+1})\ar[l]}
\]
Since $\cG^{q+1}=0$, we deduce, by descending induction on $i$,
that for all $i>0$,
\begin{equation}\label{h1dies}
H^1_{cdh}(X,\cG^i)=0
\end{equation}
Consider the diagram
\[
\xymatrix{0\ar[r]&\cG^{i+1}(R)\ar[r]\ar[d]&\cG^{i}(R)\ar[d]\ar[r]&
  \Omega^i_{F/k}\otimes\Omega^{q-i}_{R/F}\ar[r]\ar[d]&0\\
0\ar[r]&H^0_{cdh}(X,\cG^{i+1})\ar[r]&H^0_{cdh}(X,\cG^{i})\ar[r]&
  \Omega^i_{F/k}\otimes H^0_{cdh}(X,\Omega^{q-i}_{/F})\ar[r]&0}
\]
Using descending induction on $i$ again, we obtain from this diagram that
\begin{equation}\label{giagree}
\cG^i(R)\cong H^0_{cdh}(X,\cG^i) \qquad (i>0).
\end{equation}
We have a map of exact sequences
\[
\xymatrix{0\ar[r]& {}^{}\cG^1(R)\ar[d]\ar[r]&\Omega^q_{R/k}\ar[d]\ar[r]&
  \Omega^q_{R/F}\ar[r]\ar[d]&0\\
          0\ar[r]&H^0_{cdh}(X,\cG^1)\ar[r]&H^0_{cdh}(X,\Omega^q_{/k})\ar[r]&
  H^0_{cdh}(X,\Omega^q_{/F})\ar[r]&0}
\]
The third map in the bottom row is onto by \eqref{h1dies}. The
first vertical map is an isomorphism by \eqref{giagree}.  By
\eqref{preagree}, the second is an isomorphism if $q\le n$ and
onto if $q=n+1$. It follows that the same is true of the third
vertical map.
\end{proof}

\begin{prop}\label{prop:main2}
Assume $n\geq 0$. If $HH(/k)$ satisfies $n$-$cdh$-descent on $R$,
then so does $HH(/F)$.
\end{prop}

\begin{proof}
By Lemma \ref{lem:tech1}, the hypothesis is equivalent to saying
that the following conditions hold simultaneously.
\addtocounter{equation}{-1}
\begin{subequations}
\begin{gather}
HH_m^{(q)}(R/k)=0 \text{\quad if } 0\le q<m\le n\label{hhvanish}\\
\Omega^q_{R/k}\rightarrow H^0_{cdh}(X,\Omega^q_{/k})\quad
\text{is bijective if } q\le n \text{ and onto if } q=n+1\label{agree}\\
H^p_{cdh}(X,\Omega^q_{/k})=0 \quad\text{ if } p>0 \text{ and }
q\le p+n+1.\label{cdhvanish}
\end{gather}
\end{subequations}
We have to prove that the following conditions hold
\begin{subequations}
\begin{gather}
HH_m^{(q)}(R/F)=0 \text{\quad if } 0\le q<m\le n\label{hhFvanish}\\
\Omega^q_{R/F}\rightarrow H^0_{cdh}(X,\Omega^q_{/F})\quad
\text{is bijective if } q\le n \text{ and onto if } q=n+1\label{agreeF}\\
H^p_{cdh}(X,\Omega^q_{/F})=0 \quad\text{ if } p>0 \text{ and }
q\le p+n+1.\label{cdhFvanish}
\end{gather}
\end{subequations}
Using \eqref{cdhvanish}, the spectral sequence of
\ref{lem:main12}, and induction, we obtain \eqref{cdhFvanish}.
Hence \eqref{agree} implies \eqref{agreeF}, by Lemma
\ref{lem:main13}. To prove \eqref{hhFvanish} we proceed by
induction on $q$. The case $q=0$ is just the fact that
$HH^{(0)}_m(A/k)=0$ for any $m>0$, any field $k$ and any
$k$-algebra $A$.  Assume $n\ge q\ge 1$, and that we have
$HH^{(q')}_m(R/F)=0$ for all $m\le n$ and all $q'<\min\{m,q\}$. By
\eqref{agreeF} and \eqref{preinfty}, the spectral sequence of
Lemma \ref{lem:main11} collapses for $j=0$ to yield
\begin{equation*}\label{einfty}
{}_qE_{i,-i}^\infty={}_qE_{i,-i}^1.
\end{equation*}
Given this, \eqref{hhFvanish} follows from \eqref{hhvanish}
by induction.
\end{proof}

\begin{lem}\label{lem:eft=ft}
Let $F$ be a field, and $R$ a local $F$-algebra essentially of finite type. Then there
exists a field $F\subset E\subset R$ such that $R$ is isomorphic to the localization of a finite type
algebra over $E$ at a maximal ideal.
\end{lem}

\begin{proof}
The hypothesis on $R$ means that there exist an $F$-algebra $A$ of finite type and a prime ideal $P\subset A$ such that $R=A_P$.
Suppose that $\dim(A) = k+h$ and $ht(P)=h$.  By Noether normalization, there is a polynomial
subring $S = F[t_1,\dotsc,t_k]$ of $A$ meeting $P$ in $0$, and the
field $R/P=A_P/PA_P$ is a finite extension of $E=F(t_1,\dotsc,t_k)$.
There is an evident inclusion of $E$ in $R$, and $R$ is the
localization of the finite type $E$-algebra $A\otimes_{S} E$ at a
maximal ideal.
\end{proof}

The results of this section, together with Theorem \ref{thm:HHcriterion},
allow us to prove the following global regularity criterion.

\begin{thm}\label{thm:HHcodim-h}
Let $k\subseteq F$ be fields of characteristic $0$, and let $R$ be
an $F$-algebra essentially of finite type.

\noindent If $HH(/k)$ satisfies $h$-$cdh$-descent on $R$, then $R$
is smooth in codimension $h$. That is, for every prime ideal $P$
of height $h$, the local ring $R_P$ is regular.
\end{thm}

\begin{proof}
Let $P$ be a prime ideal of $R$ of height $h$. Since $HH(/F)$
satisfies $h$-$cdh$ descent on $R$ by Proposition
\ref{prop:main2}, it also satisfies $h$-$cdh$ descent on the
localization $R_P$, by Corollary \ref{cor:qc}. By Lemma
\ref{lem:eft=ft}, there is a field $F\subset E\subset R_P$ such
that $R_P$ is the localization at a maximal ideal of an algebra of
finite type over $E$. By \ref{prop:main2}, $HH(/E)$ satisfies
$h$-$cdh$ descent on $R_P$. Because $\dim (R_P)=h$, Theorem
\ref{thm:HHcriterion} implies that $R_P$ is smooth over $E$, and
hence is regular.
\end{proof}

Once again, let $F$ be a field of characteristic $0$.

\begin{thm}\label{thm:main}
Suppose $R$ is an $F$-algebra essentially of finite type. If $R$
is $K_{h+1}$-regular for some $h\geq 0$, then $R$ is regular in
codimension $h$. In particular, if $R$ is
$K_{\mathrm{dim}(R)+1}$-regular, then $R$ is regular, and hence
smooth over $F$.
\end{thm}

\begin{proof} If $R$ is $K_{h+1}$-regular, then $HH(/\Q)$ satisfies
$h$-descent on $R$, by Corollary \ref{lem:main1}. The assertion now
follows from Theorem \ref{thm:HHcodim-h}.
\end{proof}

\section{A nonreduced scheme which is $K$-regular}\label{sec:counteregg}

This section is devoted to the counterexample stated in Theorem
\ref{thm:counteregg-intro}, which reappears here as Theorem
\ref{thm:counteregg}.

\begin{lem}\label{lem:counteregg1}
Let $X$ be a smooth projective elliptic curve over a field $F$
with basepoint $Q$, and let $L$ be a degree zero line bundle $L$ on $X$.
Assume that $L$ is not an element of odd order in the Picard group.
Then $H^*(X,L^{\otimes 2n+1})=0$ for all $n\in\Z$.
\end{lem}
\begin{proof}
Because $\cO_X$ is a dualizing sheaf, we are reduced by Serre
duality to proving that $H^0(X,L^{2n+1})=0$. Because $X$ is elliptic,
there exists a rational point $P\in X$ such that $L:=L(P-Q)$.  Now if
$H^0(X,L^{\otimes2n+1})$ were nonzero, there would exist an element $f$
in the function field of $X$ with ${\rm div}(f)=(2n+1)(P-Q)$.
But because $L$ is not an odd torsion element, there is no such
$f$.
\end{proof}

\begin{thm}\label{thm:counteregg}
Let $X$ be a smooth projective elliptic curve over a field $F$ of
characteristic~0, and let $L$ be as in Lemma
\ref{lem:counteregg1}. Write $Y$ for the nonreduced scheme with
the same underlying space as $X$ but with structure sheaf
$\cO_Y=\cO_X\oplus L$, where $L$ is regarded as a square-zero
ideal.

Then for all $n$, $K_n(Y)=K_n(X)$ and $Y$ is $K_n$-regular.
\end{thm}

\begin{proof} As $X$ is regular, and hence $K_n$-regular, it suffices
to show that $K(Y\times\A^m)\to K(X\times\A^m) \cong K(X)$ is an
equivalence for all $m\ge0$. We shall prove the equivalent
assertion that the relative homotopy groups $K(Y\times\A^m,
X\times\A^m)$ are zero. By Goodwillie's theorem \cite{Goodw} and
Zariski descent, these relative $K$-groups are isomorphic to the
corresponding relative cyclic homology groups over $\Q$. By
base-change (see \cite{Kassel}) it suffices to show that the
relative groups $HH_n^{rel}=HH_n(Y,X)$ vanish for all $n$.  By
Zariski descent, it suffices to show that $H^0(X,HH_n^{rel})$ and
$H^1(X,HH_n^{rel})$ vanish for all $n$. From Lemma
\ref{lem:counteregg2} below and the fact that
$\Omega^1_{X/F}\cong\cO_X$, we see that the Zariski sheaves
$HH_n^{rel}$ are sums of odd tensor powers of $L$ when $F$ is a
number field, and odd tensor powers of $L$ tensored over $F$ with
vector spaces $\Omega^i_{F/\Q}$ in general. But the cohomology of
such powers vanishes by Lemma \ref{lem:counteregg1}.
\end{proof}

The following lemma is well-known, at least in the case when $L$ is free.
We include a proof for the sake of completeness.
For simplicity, we write $HH_*(R)$ for $HH_*(R/k)$.

\begin{lem}\label{lem:counteregg2}
Let $k$ be a field with $char(k)\ne 2$, $R$ a commutative
$k$-algebra, and $L$ a projective $R$-module of rank $1$. Put
$A=R\oplus L$. Let $M_*$ denote the graded $R$-module
\begin{equation*} M_p =
\begin{cases}L^{\otimes_R (p+1)} & p\ge0 \text{ even,}\\
L^{\otimes_R p} & p>0 \text{ odd.}\end{cases}
\end{equation*}
Then for relative Hochschild homology over $k$,
\[ HH_n(A,L)= \bigoplus_{p+q=n} M_p \otimes_R HH_q(R).\]
\end{lem}

\begin{proof}
Let $C_*(A,L)$ be the relative Hochschild complex; the subspace
$L^{\otimes_k2m+1}$ of $C_{2m}(A,L)$ consists of cycles, and
induces a map $M_{2m}=L^{\otimes_R 2m+1}\to HH_{2m}(A,L)$, because
for $x_i\in L$ and $r\in R$ we have
\[ (-1)^i b(x_0\otimes\cdots x_i\otimes r \otimes x_{i+1}\cdots) =
(x_0\otimes\cdots x_ir\otimes x_{i+1}\cdots) -
(x_0\otimes\cdots x_i\otimes rx_{i+1}\cdots).
\]
Because $Bb+bB=0$, where $B:C_{*}(A,L)\to C_{*+1}(A,L)$ is the
Connes operator, the subspace
$B(L^{\otimes_k2m+1})$ of $C_{2m+1}(A,L)$ also consists of cycles
and induces a map $M_{2m+1}=L^{\otimes_R 2m+1}\to HH_{2m}(A,L)$.
Thus we have a graded map $M_*\to HH_{2m}(A,L)$. Because $HH_*(A,L)$
is a graded module over $HH_*(R)$, we get a canonical $R$-module map from
$M_*\otimes_R HH_*(R)$ to $HH_*(A,L)$. To see that it is an isomorphism,
we may assume $R$ is local, whence $A=R[x]/<x^2>$. By the K\"unneth
formula, we are reduced to the case $R=k$, which is straightforward.
\end{proof}

\subsection*{Acknowledgements}
The authors would like to thank M.~Schlichting, whose contributions
go beyond the collaboration \cite{chsw}.
C.~ Weibel would like to thank his children for riding the roller coaster
``Dueling Dragons" enough times for him to work out the example
in Theorem \ref{thm:counteregg-intro}.

\end{document}